\documentclass[a4paper, 12pt]{amsart}
\usepackage{amssymb, amstext, amscd, amsmath, color}

\begin{document}

\newtheorem{theorem}{Theorem}[section]
\newtheorem{corollary}[theorem]{Corollary}
\newtheorem{proposition}[theorem]{Proposition}
\newtheorem{lemma}[theorem]{Lemma}
%
%      Theorem style with roman text, numbered within section

\theoremstyle{definition}
\newtheorem{remark}[theorem]{Remark}
\newtheorem{definition}[theorem]{Definition}
\newtheorem{note}[theorem]{Note}
\newtheorem{example}[theorem]{Example}
\newcommand{\Prf}{\noindent\textbf{Proof.\ }}
\newcommand{\bx}{\hfill$\blacksquare$\medbreak}
\newcommand{\upbx}{\vspace{-2.5\baselineskip}\newline\hbox{}%
\hfill$\blacksquare$\newline\medbreak}
\newcommand{\eqbx}[1]{\medbreak\hfill\(\displaystyle #1\)\bx}
% for proofs ending in a one line equation
%      Useful shortforms

\newcommand{\FFock}{\mathcal{F}}
\newcommand{\kil}{\mathsf{k}}
\newcommand{\Hil}{\mathsf{H}}
\newcommand{\hil}{\mathsf{h}}
\newcommand{\Kil}{\mathsf{K}}
\newcommand{\Real}{\mathbb{R}}
\newcommand{\Rplus}{\Real_+}

%      Blackboard bold letters
\newcommand{\bC}{{\mathbb{C}}}
\newcommand{\bD}{{\mathbb{D}}}
\newcommand{\bN}{{\mathbb{N}}}
\newcommand{\bQ}{{\mathbb{Q}}}
\newcommand{\bR}{{\mathbb{R}}}
\newcommand{\bT}{{\mathbb{T}}}
\newcommand{\bX}{{\mathbb{X}}}
\newcommand{\bZ}{{\mathbb{Z}}}
\newcommand{\bH}{{\mathbb{H}}}
%      Useful shortforms
\newcommand{\BH}{{\B(\H)}}
\newcommand{\bsl}{\setminus}
\newcommand{\ca}{\mathrm{C}^*}
\newcommand{\cstar}{\mathrm{C}^*}
\newcommand{\cenv}{\mathrm{C}^*_{\text{env}}}
\newcommand{\rip}{\rangle}
\newcommand{\ol}{\overline}
\newcommand{\td}{\widetilde}
\newcommand{\wh}{\widehat}
\newcommand{\sot}{\textsc{sot}}
\newcommand{\wot}{\textsc{wot}}
\newcommand{\wotclos}[1]{\ol{#1}^{\textsc{wot}}}
%      Capital script letters
 \newcommand{\A}{{\mathcal{A}}}
 \newcommand{\B}{{\mathcal{B}}}
 \newcommand{\C}{{\mathcal{C}}}
 \newcommand{\D}{{\mathcal{D}}}
 \newcommand{\E}{{\mathcal{E}}}
 \newcommand{\F}{{\mathcal{F}}}
 \newcommand{\G}{{\mathcal{G}}}
\renewcommand{\H}{{\mathcal{H}}}
 \newcommand{\I}{{\mathcal{I}}}
 \newcommand{\J}{{\mathcal{J}}}
 \newcommand{\K}{{\mathcal{K}}}
\renewcommand{\L}{{\mathcal{L}}}
 \newcommand{\M}{{\mathcal{M}}}
 \newcommand{\N}{{\mathcal{N}}}
\renewcommand{\O}{{\mathcal{O}}}
\renewcommand{\P}{{\mathcal{P}}}
 \newcommand{\Q}{{\mathcal{Q}}}
 \newcommand{\R}{{\mathcal{R}}}
\renewcommand{\S}{{\mathcal{S}}}
 \newcommand{\T}{{\mathcal{T}}}
 \newcommand{\U}{{\mathcal{U}}}
 \newcommand{\V}{{\mathcal{V}}}
 \newcommand{\W}{{\mathcal{W}}}
 \newcommand{\X}{{\mathcal{X}}}
 \newcommand{\Y}{{\mathcal{Y}}}
 \newcommand{\Z}{{\mathcal{Z}}}

\newcommand{\sgn}{\operatorname{sgn}}
\newcommand{\rank}{\operatorname{rank}}
\newcommand{\Isom}{\operatorname{Isom}}
\newcommand{\qIsom}{\operatorname{q-Isom}}
\newcommand{\Cknet}{{\mathcal{C}_{\text{knet}}}}
\newcommand{\Ckag}{{\mathcal{C}_{\text{kag}}}}
\newcommand{\rind}{\operatorname{r-ind}}
\newcommand{\lind}{\operatorname{r-ind}}

% % TOP MATTER

\title[Infinitesimal rigidity for non-Euclidean frameworks]% end with percent
 {Infinitesimal rigidity for non-Euclidean bar-joint frameworks} % This is the full title of the paper
% Use lowercase letters in title except for proper names
% Avoid equations in title if possible
% Do not use the \thanks{} command; use \extraline{} instead (see below).

\author[D. Kitson and S. C. Power]{D. Kitson and S. C. Power}
\thanks{Supported by EPSRC grant  EP/J008648/1.}
\address{Dept.\ Math.\ Stats.\\ Lancaster University\\
Lancaster LA1 4YF \\U.K. }
\email{d.kitson@lancaster.ac.uk, s.power@lancaster.ac.uk}
%\dedication{A dedication can be included here}

%Insert `2000 Mathematics Subject Classification' numbers here:

\begin{abstract}
The minimal infinitesimal rigidity of  bar-joint frameworks in the non-Euclidean spaces $(\mathbb{R}^2, \|\cdot \|_q)$ for $1\leq q \leq \infty, q \neq 2$, are characterised in terms of $(2,2)$-tight graphs.
Specifically, a generically placed bar-joint framework $(G,p)$ in the plane is minimally infinitesimally rigid with respect to a non-Euclidean $\ell^q$ norm if and only if the underlying graph $G=(V,E)$ contains $2|V|-2$ edges and every subgraph $H=(V(H),E(H))$ contains at most $2|V(H)|-2$ edges.

\end{abstract}

\maketitle

%%%%%%%%%%%%%%%%%%%%%%%%%%%%%%%%%%%%%%%%%%%%%%%%%%%%%%%%%%%%%%
\section{Introduction}
It is a longstanding open problem in the  theory of bar-joint frameworks
to obtain a form of combinatorial characterisation for the infinitesimal rigidity of generic frameworks in the Euclidean space $\mathbb{R}^3$. 
On the other hand in two dimensions the foundational characterisation of Laman \cite{Lam} provides a necessary and sufficient condition,
namely that the underlying simple graph $(V,E)$  should contain a spanning subgraph $(V,E')$ which is $(2,3)$-tight, meaning that the Maxwell count
$|E'|=2|V|-3$ should hold, together with the inequalities $|E(H)|\leq 2|V(H)|-3$ for all subgraphs $H$ with $|E(H)|\geq 1$. 
We show that for the  non-Euclidean norms $\|\cdot\|_q, 1 \leq q \leq \infty, q\neq 2,$ 
there are exact analogues of Laman's characterisation in terms of $(2,2)$-tight graphs and in this spirit we pose several problems on the characterisation of generic bar-joint frameworks in normed spaces.

For the  theory of bar-joint frameworks in finite dimensional Euclidean space see, for example, Asimow and Roth \cite{A&R}, \cite{A&R_1979}, Gluck \cite{Gluck} and Graver, Servatius and Servatius \cite{GSS}. 
It seems that the setting of general finite-dimensional normed spaces is new and requires
novel combinatorial graph theory. 
One can verify for example that a triangular framework $(K_3, p)$ in $(\mathbb{R}^2, \|\cdot\|_q)$
has a flex (of rotational type) which does not derive from an ambient rigid motion and
that a regular $K_4$ framework is minimally infinitesimally rigid. 
These observations echo similar facts for bar-joint frameworks in three dimensions whose vertices are constrained
to a circular cylinder, and indeed the circular cylinder and $(\mathbb{R}^2, \|\cdot\|_q)$ both possess a two
dimensional space of trivial infinitesimal motions. Accordingly it is $K_4$ (or $K_1$) rather than
$K_2$ which plays the role of a base graph in the inductive scheme appropriate for these
normed spaces.

A standard proof of Laman's theorem makes use of the two Henneberg moves $G \to G'$ which increase by one the number of vertices and which generate all $(2,3)$-tight graphs from the base graph $K_2$. For $q\not=1,2,\infty$ we take a similar approach using both an associated rigidity matrix and an inductive scheme for $(2,2)$-tight graphs in terms of the two Henneberg moves and two additional moves, namely the vertex-to-$4$-cycle move and the vertex-to-$K_4$ move. 
For $q\in\{1,\infty\}$ some novel features emerge in that infinitesimal rigidity is characterised by induced graph colourings and the set of regular realisations is no longer dense. We develop this idea in the broader context of polytopic norms.

%%%%%%%%%%%%%%%%%%%%%%%%%%%%%%%%%%%%%%%%%%%%%%%%%%%%%%%%%%%%%%%
\section{Infinitesimal flexes}
We define a bar-joint framework in the normed linear space $(X, \|\cdot\|)$ to be a pair $(G, p)$ with $G$
a simple graph $(V,E)$ and $p = (p_1,p_2,\ldots, p_n)$ an $n$-tuple of points in $X$ representing the
placement of the vertices of $V$ with respect to some given labelling $v_1,v_2,\ldots, v_n$ of $V$. 
In fact we may assume that the graph is connected and so we do so throughout.

\begin{definition}
An {\em infinitesimal flex} of $(G, p)$ is a vector $u = (u_1, \ldots, u_n)\in X^n$  such that the corresponding framework edge lengths of $(G, p+tu)$ have second order deviation
from the original lengths,
\[\|(p_i+tu_i)-(p_j+tu_j)\|-\|p_i-p_j\|=o(t), \,\,\,\,\,\,\,\, \mbox{ as } t\to 0\]
for each edge $v_iv_j\in E$.
\end{definition}

The infinitesimal flexes of a framework $(G,p)$ form a proper closed linear subspace of $X^n$.

\begin{definition}
A {\em rigid motion} of the normed space $(X,\|\cdot\|)$ is a collection $\{\gamma_x\}_{x\in X}$ of continuous paths $\gamma_x:[-1,1]\to X$  such that
\begin{enumerate}
\item
$\gamma_x$ is differentiable at $0$ and
$\gamma_x(0)=x$ for each  $x\in X$, and,
\item
$\|\gamma_x(t)-\gamma_y(t)\| = \|x-y\|$ for all $t\in [-1,1]$ and for all $x,y\in X$.
\end{enumerate}
\end{definition}

Each rigid motion gives rise to a vector field on $X$ via the map $x\mapsto \gamma_x'(0)$ and this vector field induces an infinitesimal flex on any given framework $(G,p)$ as shown in the following lemma. An infinitesimal flex is said to be {\em trivial} if it is obtained from a rigid motion  in this way.

\begin{lemma}
\label{TrivialFlex}
Let $(G,p)$ be a bar-joint framework in a  normed linear space $(X,\|\cdot\|)$.
\begin{enumerate}
\item
If $\{\gamma_x\}_{x\in X}$ is a rigid motion then $u=(\gamma'_{p_1}(0),\ldots,\gamma'_{p_n}(0))$ is a (trivial) infinitesimal flex of $(G,p)$.
\item
If $a\in X$ then $u=(a,\ldots,a)\in X^{n}$ is a trivial infinitesimal flex of $(G,p)$.
\end{enumerate}
\end{lemma}

\proof
$(i)$
For each edge $v_iv_j\in E$ we compute
\begin{eqnarray*}
\epsilon(t) &=&\left| \frac{\|(p_i+tu_i)-(p_j+tu_j)\|-\|p_i-p_j\|}{t} \right| \\ 
&=&  \left| \frac{\|(\gamma_{p_i}(0)+t\gamma_{p_i}'(0))-(\gamma_{p_j}(0)+t\gamma_{p_j}'(0))\|-\|\gamma_{p_i}(t)-\gamma_{p_j}(t)\|}{t} \right|\\ 
&\leq& 
\left\| \frac{(\gamma_{p_i}(t)-\gamma_{p_j}(t))-(\gamma_{p_i}(0)-\gamma_{p_j}(0))}{t}-(\gamma_{p_i}'(0)-\gamma_{p_j}'(0))\right\| 
\end{eqnarray*}
From the latter quantity we see that $\epsilon(t)\to 0$ as $t\to 0 $ and so $u$ is an infinitesimal flex.

$(ii)$
Clearly $u=(a,\ldots,a)$ is an infinitesimal flex and is induced by 
the rigid motion $\{\gamma_x\}_{x\in X}$ with $\gamma_x(t)=x+at$ for each $x\in X$.
\endproof

We recall the following theorem of Mazur and Ulam.

\begin{theorem}[(\cite{MU})]
If $A:X\to Y$ is a surjective isometry between real normed linear spaces $X$ and $Y$ and $A(0)=0$
then $A$ is a linear map.
\end{theorem}

In the case of the $\ell^q$ norms and polytopic norms on $\mathbb{R}^d$ considered in Sections \ref{lqNorms} and \ref{PolytopicNorms} the trivial infinitesimal flexes of a bar-joint framework are described by the following lemma.

\begin{lemma}
\label{RigidMotions}
Let $(G,p)$ be a bar-joint framework in $(X,\|\cdot\|)$ where $X$ is a finite dimensional normed linear space over $\mathbb{R}$ which admits only finitely many surjective linear isometries.
Then $u$ is a trivial infinitesimal flex of $(G,p)$  if and only if $u=(a,\ldots,a)\in X^{n}$ for some $a\in X$.
\end{lemma}

\proof
The sufficiency of the condition is proved in Lemma \ref{TrivialFlex}$(ii)$. To prove necessity 
suppose that $u$ is a trivial infinitesimal flex of $(G,p)$ and let $\{\gamma_x\}_{x\in X}$ be a rigid motion with $u=(\gamma_{p_1}'(0),\ldots,\gamma_{p_n}'(0))$.
By definition the mapping $\Gamma_t:X\to X$, $x\mapsto \gamma_x(t)$ is an isometry for each $t\in [-1,1]$.
Moreover, since $X$ is finite dimensional each $\Gamma_t$ is a surjective isometry.
Define $A_t:X\to X$, $x\mapsto \Gamma_t(x)-\Gamma_t(0)$ for each $t\in [-1,1]$. 
Then $A_t$ is a surjective isometry with $A(0)=0$ and so, by the Mazur-Ulam theorem, $A_t$ is linear.
Let $I,T_1,\ldots,T_m$ be the finitely many surjective linear isometries on $X$ and choose vectors
 $x_1,\ldots,x_m\in X$ with $T_j(x_j)\not=x_j$. 
Since $\gamma_x$ is continuous it follows that $A_t(x)\to x$ as $t\to0$ for all $x\in X$. 
Thus if  we set $\epsilon =  \min_{j=1,\ldots,m}\|T_j(x_j)-x_j\|$ then
there exists $\delta>0$ such that $\max_{j=1,\ldots,m} \|A_t(x_j)-x_j\|<\epsilon$ for all $|t|<\delta$.
We conclude that $A_t=I$ for all $|t|<\delta$.
From the definition of $A_t$ we now have  $\Gamma_t(x)=x+\Gamma_t(0)$ for all $|t|<\delta$ and so 
$\gamma'_x(0)=\gamma_0'(0)$ for all $x\in X$.
In particular,  $u=(\gamma_0'(0),\ldots,\gamma_0'(0))\in X^n$.

\endproof

\begin{definition}
A bar-joint framework $(G,p)$ is {\em infinitesimally flexible} in $(X,\|\cdot\|)$ if it has a non-trivial infinitesimal flex. Otherwise it is said to be {\em infinitesimally rigid}.
\end{definition}

A framework $(G,p)$ is {\em minimally infinitesimally rigid} (or {\em isostatic}) if it is infinitesimally rigid and removing a single edge from $G$ results in an infinitesimally flexible framework.

\begin{lemma}
\label{Isometric}
Let $A:X\to Y$ be an isometric affine isomorphism between normed linear spaces $(X,\|\cdot\|_X)$ and $(Y,\|\cdot\|_Y)$.
Then a bar-joint framework $(G,p)$  in $X$ is (minimally) infinitesimally  rigid if and only if $(G,A(p))$ is (minimally) infinitesimally  rigid in $Y$.
\end{lemma}

\proof
Here we are using the obvious notation $A(p)=(A(p_1),\ldots,A(p_n))$.
If a vector $u\in X^n$ is an infinitesimal flex for $(G,p)$  then $A(u)\in Y^n$ is clearly an infinitesimal flex for $(G,A(p))$ and vice versa.
There is a one-to-one correspondence between the rigid motions $\{\gamma_x\}_{x\in X}$ of $X$ 
and the rigid motions $\{\tilde{\gamma}_y\}_{y\in Y}$ of $Y$ given by 
$\tilde{\gamma}_{A(x)}(t) = A(\gamma_{x}(t))$.
If $u=(\gamma_{p_1}'(0),\ldots,\gamma_{p_n}'(0))$ is a trivial infinitesimal flex  of $(G,p)$ then $A(u)=(\tilde{\gamma}_{A(p_1)}'(0),\ldots, \tilde{\gamma}_{A(p_n)}'(0))$ and so $A(u)$ is a trivial infinitesimal flex of $(G,A(p))$.
We conclude that if $(G,p)$ is infinitesimally rigid in $X$ then $(G,A(p))$ is infinitesimally rigid in $Y$.
\endproof

%%%%%%%%%%%%%%%%%%%%%%%%%%%%%%%%%%%%%%%%%%%%%%%%%%%%%%%%%%%%%%%%
\section{$\ell^q$ norms}
\label{lqNorms}
The aim of this section is to introduce the rigidity matrix for a framework $(G,p)$  in $(\mathbb{R}^d,\|\cdot\|_q)$  and to characterise the minimally infinitesimally rigid frameworks in $(\mathbb{R}^2,\|\cdot\|_q)$ when $1<q<\infty$ and $q\not=2$.
The $1$-norm and $\infty$-norm are dealt with in the next section.

For a point $a=(a_1,\ldots,a_d)\in\mathbb{R}^d$ and $k\in(0,\infty)$ we write 
\[a^{(k)}=(sgn(a_1)|a_1|^k,\ldots,sgn(a_d)|a_d|^k)\]
where $sgn$ is the sign function.
We remark that for $1<q<\infty$ the mapping $\mathbb{R}^d\to\mathbb{R}^d$, $a\mapsto  a^{(q-1)}$ is injective. It follows that  $p_1,p_2,p_3$ are non-collinear points in $\mathbb{R}^d$ if and only if $(p_1-p_2)^{(q-1)}$ and $(p_1-p_3)^{(q-1)}$ are linearly independent. We will make use of this fact later.

\begin{definition}
The {\em rigidity matrix} $R_q(G,p)$ is a $|E|\times nd$ matrix with rows indexed by the edges of $G$ and $nd$ many columns indexed by the coordinates of the vertex placements $p_1,\ldots,p_n$. The row entries which correspond to an edge $v_iv_j\in E$ are
\[
\left[\begin{array}{ccccccccccc}
 0 & \cdots& 0& (p_i-p_j)^{(q-1)} & 0 &\cdots & 0 & -(p_i-p_j)^{(q-1)} & 0 & \cdots & 0
 \end{array}\right]
 \]
 where  non-zero entries may only appear in the $p_i$ columns and the $p_j$ columns.
\end{definition}

\begin{proposition}
\label{RigidityMatrix}
Let $(G,p)$ be a bar-joint framework in $(\mathbb{R}^d,\|\cdot\|_q)$, $1< q <\infty$, $q\not=2$.
Then 
\begin{enumerate}
\item
$u\in \mathbb{R}^{nd}$ is an infinitesimal flex for $(G,p)$  if and only if $R_q(G,p)u=0$.
\item
$(G,p)$  is infinitesimally rigid if and only if 
$\rank R_q(G,p) = dn-d$.
\end{enumerate}
\end{proposition}

\proof 
$(i)$
Let $u=(u_1,\ldots,u_n)\in \mathbb{R}^{nd}$ and for each edge $v_iv_j\in E$   consider the associated edge length function
\[\zeta_{ij}:\mathbb{R}\to\mathbb{R}, \,\,\,\,\,\, t\mapsto \|(p_i+tu_i)-(p_j+tu_j)\|_q\]
Note that $u$ is an infinitesimal flex for $(G,p)$ if and only if $\zeta_{ij}'(0)=0$ for each edge $v_iv_j\in E$. 
We  compute
\[\zeta_{ij}'(0)=\frac{(p_i-p_j)^{(q-1)}\cdot(u_i-u_j)}{\|p_i-p_j\|_q^{q-1}}\]
from which the result follows.

$(ii)$
 Recall that there are only finitely many surjective linear isometries on $(\mathbb{R}^d,\|\cdot\|_q)$ for $q\in(1,\infty)$, $q\not=2$. These are  given by signed permutation matrices (see \cite[Proposition 2.f.14]{LT}) and  have the form $T(x)=(\theta_1x_{\pi(1)},\ldots,\theta_d x_{\pi(d)})$ for some permutation $\pi$
of $\{1,2,\ldots,d\}$ and some $\theta_j\in\{-1,1\}$. Applying Lemma \ref{TrivialFlex} we see that the trivial infinitesimal flexes form a $d$ dimensional subspace of $\mathbb{R}^{nd}$. The result now follows from $(i)$.
\endproof

As in Euclidean space the set of all realisations of a bar-joint framework in $(\mathbb{R}^d,\|\cdot \|_q)$ can be partitioned into regular and non-regular realisations. 

\begin{definition}
A framework $(G,p)$ is \emph{regular} in $(\mathbb{R}^d,\|\cdot \|_q)$ if the rank of its rigidity matrix $R_q(G,p)$ is maximal over all framework realisations $(G,p')$. 
\end{definition}

Note that for $q\in(1,\infty)$ the set of all regular realisations for $G$ is an open dense subset of $\mathbb{R}^{nd}$.
If a bar-joint framework $(G,p)$ is infinitesimally rigid in  $(\mathbb{R}^d,\|\cdot\|_q)$ then $p$ is necessarily a regular realisation.
If a framework $(G,p)$ is infinitesimally flexible for all realisations $p$ in some open subset of $\mathbb{R}^d$ then $(G,p)$ must be infinitesimally flexible for all realisations in $\mathbb{R}^d$.

\begin{example}
It follows from  Proposition \ref{RigidityMatrix} that any realisation of the complete graphs $K_2$ and $K_3$
in $(\mathbb{R}^2,\|\cdot\|_q)$  with $1< q<\infty$, $q\not=2$, must be infinitesimally 
flexible. In contrast by
computing an appropriate rigidity matrix we see that any regular realisation of $K_4$ in
these spaces will be infinitesimally rigid and indeed minimally infinitesimally rigid. 
 With our definitions note that any single vertex framework $(K_1,p)$ is infinitesimally rigid.
\end{example}

\begin{remark}
That a regular $K_3$ framework is infinitesimally  flexible may seem something of a paradox
and in fact there is a measure of freedom in how one may view the flexibility of small
complete graphs in finite-dimensional normed spaces. One could adopt a more intrinsic
notion of infinitesimal rigidity (rather than the ambient isometry one we have given above)
and assert that a framework $(G, p)$ is infinitesimally rigid if there is no infinitesimal 
flex $(u_1,\ldots, u_n)$ for which $\zeta'_{ij}(0)\not=0$ for some non-edge pair $v_i, v_j$. By this requirement a regular $K_3$
framework is infinitesimally rigid $(\mathbb{R}^2,\|\cdot\|_q)$ for vacuous reasons. On the other hand the double triangle resulting from a Henneberg-$1$ move (defined below) applied to $K_3$ is infinitesimally flexible by this definition (as well as the one above). As a Henneberg $1$-move in general preserves infinitesimal rigidity, $K_3$ is exceptional in this sense also.
\end{remark}

A simple connected graph $G=(V,E)$ is said to be {\em $(2,2)$-sparse} if 
$|E(H)|\leq 2|V(H)|-2$ for all subgraphs $H=(V(H),E(H))$.
If in addition $|E|=2|V|-2$ then $G$ is said to be {\em $(2,2)$-tight}.
See also Lee and Streinu \cite{Lee&Streinu} and Szeg\"{o} \cite{Szego}.

We now state the main result of this section (Theorem \ref{qNormThm}) which is a characterisation of the minimally infinitesimally rigid bar-joint frameworks
in $(\mathbb{R}^2,\|\cdot\|_q)$ for $q\in (1,\infty)$ and $q\not=2$.

\begin{theorem}
\label{qNormThm}
Let $(G,p)$ be a  regular bar-joint framework in $(\mathbb{R}^2,\|\cdot\|_q)$, $1< q<\infty$, $q\not=2$.
Then $(G,p)$ is minimally infinitesimally  rigid if and only if $G$  is $(2,2)$-tight.
\end{theorem}

The rigidity matrix will play a prominent role in the proof of Theorem \ref{qNormThm} and necessity can be proved directly. To prove sufficiency we will use an  inductive construction based on certain graph moves.

\begin{definition}
Let $G=(V,E)$ be a simple connected graph. A graph $G'=(V',E')$ is obtained from $G$ by a {\em Henneberg $1$-move} by the following process: 
\begin{enumerate}
\item Adjoin a new vertex $v_{0}$ to $V$ so that $V'=V\cup\{v_0\}$.
\item Choose two distinct vertices $v_1,v_2\in V$ and adjoin the edges $v_0v_{1}$ and $v_0v_{2}$ to $E$ so that
$E'=E\cup \{v_0v_1,\,v_0v_2\}$.
\end{enumerate}
\end{definition}

\begin{lemma}
\label{Henneburg1}
 Let $(G,p)$ be an infinitesimally rigid bar-joint framework in $(\mathbb{R}^2,\|\cdot\|_q)$, $1<q <\infty$, $q\not=2$. 
If $G'$ is obtained from $G$ by a Henneberg $1$-move then
 $(G',p')$ is infinitesimally  rigid  for some realisation $p'$. 
\end{lemma}

\proof
Suppose the Henneberg $1$-move $G\to G'$  is based on the vertices $v_1,v_2\in V$.
Choose $p'=(p_0,p_1,\ldots,p_n)$ such that $p'$ is a realisation of $G'$ with $p_0,p_1,p_2$ not collinear.
If $(G',p')$ has an infinitesimal flex $u' = (u_0,u_1,\ldots,u_n)$ then  $u=(u_1,\ldots,u_n)$ is an infinitesimal flex for $(G,p)$. 
Thus $u$ is trivial and by Lemma \ref{RigidMotions}, $u_1=\cdots=u_n$.
Now consider the rows of the rigidity matrix $R_q(G',p')$ which correspond to the edges $v_0v_{1}$ and $v_0v_{2}$,
\begin{eqnarray}
\label{RigEqn}
\left[\begin{array}{cccccc}
(p_0-p_1)^{(q-1)} & -(p_0-p_1)^{(q-1)} & 0 & 0 & \cdots & 0 \\
(p_0-p_2)^{(q-1)} & 0 & -(p_0-p_2)^{(q-1)} & 0 & \cdots & 0
\end{array}\right]\end{eqnarray}
Note that $u_0-u_{1}=u_0-u_{2}$ is orthogonal (in the Euclidean sense) to both  $(p_0-p_{1})^{(q-1)}$ and  $ (p_0-p_{2})^{(q-1)}$.
As remarked previously the latter vectors are linearly independent in $\mathbb{R}^2$ and so $u_{0}=u_1$.
We conclude that $u'$ is trivial and so $(G',p')$ is infinitesimally rigid.
\endproof

\begin{definition}
Let $G=(V,E)$ be a simple connected graph. A graph $G'=(V',E')$ is obtained from $G$ by a {\em Henneberg $2$-move} by the following process:  
\begin{enumerate}
\item Choose an edge $v_1v_2\in E$ and one other vertex $v_3\in V$.
\item Adjoin a new vertex $v_0$ to $V$ so that $V'=V\cup \{v_0\}$.
\item Remove the edge $v_1v_2$ from $E$ and adjoin the edges $v_0v_1$, $v_0v_2$ and $v_0v_3$ so that
$E'=(E\backslash\{v_1v_2\})\cup\{v_0v_1,\,v_0v_2,\,v_0v_3\}$.  
\end{enumerate}
\end{definition}

\begin{lemma}
\label{Henneburg2}
Let $(G,p)$ be an infinitesimally rigid bar-joint framework in $(\mathbb{R}^2,\|\cdot\|_q)$, $1<q <\infty$, $q\not=2$.
If $G'$ is obtained from $G$ by a Henneberg $2$-move 
then $(G',p')$ is infinitesimally rigid for some  realisation $p'$.
\end{lemma}

\proof
The set of infinitesimally rigid realisations of $G$  form a dense open set in $\mathbb{R}^{2n}$ and so by replacing $p=(p_1,p_2,\ldots,p_n)$ if necessary we can assume that  $p_1,p_2,p_3$ are non-collinear in $\mathbb{R}^2$.
Choose  $p'=(p_0,p_1,\ldots,p_n)$ such that $p_{0}$ is a point on the line segment joining $p_1$ to $p_2$.  
The respective rigidity matrices can be expressed in the following block form,
\[R_q(G',p') = \left[
\small{ \begin{array}{c|c}
\begin{array}{c}
 (p_0-p_1)^{(q-1)}   \\
(p_0-p_2)^{(q-1)}   \\
(p_0-p_3)^{(q-1)}  
\end{array} & \begin{array}{ccccc}
 -(p_0-p_1)^{(q-1)} & 0 & 0 & \cdots & 0  \\
 0 & -(p_0-p_2)^{(q-1)} & 0 & \cdots & 0  \\
 0 & 0 &-(p_0-p_3)^{(q-1)}   & \cdots & 0  \end{array} \\ \hline
& \\
0 &  Z \\
&
\end{array} } \right]\]
\[R_q(G,p) = \left[ \small{ \begin{array}{c}
\begin{array} {ccccc}
(p_1-p_2)^{(q-1)}  & -(p_1-p_2)^{(q-1)} & 0  & \cdots & 0   \end{array} \\ \hline
\\  Z  \\
  \end{array} }\right]
\]
The collinearity of $p_0,p_1,p_2$ and non-collinearity of $p_1,p_2,p_3$ ensure that 
\[\rank R_q(G',p')=\rank R_q(G,p)+2 = 2(n+1)-2\]
and so $(G',p')$ is infinitesimally rigid by Proposition \ref{RigidityMatrix}.
\endproof

\begin{definition}
Let $G=(V,E)$ be a simple connected graph. A graph $G'=(V',E')$ is obtained from $G$ by a {\em vertex to $4$-cycle move} by the following process:
\begin{enumerate}
\item Choose a vertex $v_1\in V$ together with two edges $v_1v_2,\, v_1v_3\in E$.
\item Adjoin a new vertex $v_0$ to $V$ so that $V'=V\cup \{v_0\}$ and adjoin the edges $v_0v_2$ and $v_0v_3$ to  $E$.
\item Either leave any remaining edge of the form $v_1w\in E$ unchanged or replace it with $v_0w$.
\end{enumerate}
\end{definition}

\begin{lemma}
\label{Vertex-to-4-cycle}
Let $(G,p)$ be an infinitesimally rigid bar-joint framework in $(\mathbb{R}^2,\|\cdot\|_q)$, $1<q <\infty$, $q\not=2$. 
If $G'$ is obtained from $G$ by a  vertex to $4$-cycle move 
then $(G',p')$ is infinitesimally  rigid  for some realisation $p'$. 
\end{lemma}

\proof
The  regular realisations of $G$ and $G'$ form dense open subsets of  $\mathbb{R}^{2n}$ and $\mathbb{R}^{2(n+1)}$ respectively.
It follows that there exists  $p'=(p_0,p_1,p_2,\ldots,p_n)\in \mathbb{R}^{2(n+1)}$ such that $p_1,p_2,p_3$ are non-collinear and $(G',p')$ and $(G,p)$ are both  regular where $p=(p_1,p_2,\ldots,p_n)$.

Since $p'$ is a regular realisation the rank of the rigidity matrix $R_q(G',p')$ does not increase when we replace $p'$ with $p^*=(p_1,p_1,p_2,\ldots,p_n)$. 
Let $u'=(u_0,u_1,\ldots,u_n)$ be in the kernel of $R_q(G',p^*)$ 
and consider the rows of  $R_q(G',p^*)$ which correspond to the edges $v_0v_2$, $v_0v_3$, $v_1v_2$ and $v_1v_3$,
\begin{eqnarray*}
\left[ \begin{array}{ccccccc}
(p_1-p_2)^{(q-1)} & 0 & -(p_1-p_2)^{(q-1)} & 0 & 0 & \cdots & 0 \\
(p_1-p_3)^{(q-1)} & 0 & 0 & -(p_1-p_3)^{(q-1)} & 0 & \cdots & 0 \\
0 & (p_1-p_2)^{(q-1)} & -(p_1-p_2)^{(q-1)} & 0 & 0 & \cdots & 0 \\
0 & (p_1-p_3)^{(q-1)} & 0  & -(p_1-p_3)^{(q-1)}  & 0 & \cdots & 0 
\end{array}\right]
\end{eqnarray*}
This sytem leads us to the orthogonality relations,
\[(p_1-p_2)^{(q-1)}\cdot(u_0-u_1)=0=(p_1-p_3)^{(q-1)}\cdot(u_0-u_1)\]
From our earlier remark, since $p_1,p_2,p_3$ are non-collinear  $(p_1-p_2)^{(q-1)}$ and $(p_1-p_3)^{(q-1)}$ are linearly independent
in $\mathbb{R}^2$ and so $u_0=u_1$.
It  follows that $R_q(G,p)(u_1,\ldots,u_n)=0$ and so $u_1=\cdots=u_n$. 
That $(G',p')$ is infinitesimally rigid is now clear from Proposition \ref{RigidityMatrix} since
\[\rank R_q(G',p^*)=2(n+1)-\dim \ker R_q(G',p^*)=2(n+1)-2 \]
\endproof

\begin{definition}
Let $G=(V,E)$ be a simple connected graph. A graph $G'=(V',E')$ is obtained from $G$ by a {\em vertex-to-$K_4$ move} by the following process:
\begin{enumerate}
\item Remove a vertex $v_1$ from $V$ and adjoin four new vertices $w_1,w_2,w_3,w_4$ so that $V'=(V\backslash\{v_1\})\cup\{w_1,w_2,w_3,w_4\}$.
\item Adjoin all edges of the form $w_iw_j$ to $E$.
\item Replace each edge of the form $v_1v\in E$ with $w_jv$ for some $j$.
\end{enumerate}
\end{definition}

\begin{lemma}
\label{Vertex-to-K4}
Let $(G,p)$ be an infinitesimally rigid bar-joint framework in $(\mathbb{R}^2,\|\cdot\|_q)$, $1<q<\infty$, $q\not=2$. 
If $G'$  is obtained from $G$ by a vertex-to-$K_4$ move 
then $(G',p')$ is infinitesimally  rigid  for some realisation $p'$. 
\end{lemma}

\proof
Let $H$ be the subgraph of $G'$ which is a copy of $K_4$ based on the new vertices $w_1,w_2,w_3,w_4$.
Given any realisation $p'=(p_H,p_{G'\backslash H})$ of $G'$ the rigidity matrix  can be expressed in block form
\begin{eqnarray}
\label{Eqn1}
R_q(G',p') = \left[ \begin{array} {cc}
R_q(H,p_H) & 0 \\
X_1(p') & X_2(p') \end{array} \right]
\end{eqnarray}
The  regular realisations of $G$,  $H$ and $G'$ form dense open sets in $\mathbb{R}^{2n}$, $\mathbb{R}^{8}$ and $\mathbb{R}^{2(n+3)}$ respectively.
It follows that there exists a realisation $p'=(p_1',p_2',p_3',p_1,p_2,\ldots,p_n)\in \mathbb{R}^{2(n+3)}$ such that the frameworks $(G,p)$, $(H,p_H)$ and $(G',p')$ are all regular and such that the matrix $X_2(p')$ achieves its maximal rank over all realisations.
Here  $p=(p_1,p_2,\ldots,p_n)$ and $p_H=(p_1',p_2',p_3',p_1)$. 

If $(G',p')$ has an infinitesimal flex $u'=(u_1',u_2',u_3',u_1,u_2,\ldots,u_n)\in\mathbb{R}^{2(n+3)}$ then by making a translation we can assume that $u_1=0$. Clearly $R_q(H,p_H)(u_1',u_2',u_3',u_1)=0$.
As was remarked earlier, every regular realisation of $K_4$ in $(\mathbb{R}^2,\|\cdot\|_q)$ is infinitesimally rigid and so using Lemma \ref{RigidMotions} we have $u_1'=u_2'=u_3'=u_1=0$.
Consider now the ill-positioned framework
$(G', p^*)$ where $p^*=(p_1,p_1,p_1,p_1,p_2,\ldots,p_n)$ and its rigidity matrix
\begin{eqnarray}
\label{Eqn1}
R_q(G',p^*) = \left[ \begin{array} {cc}
R_q(H,p'') & 0 \\
X_1(p^*) & X_2(p^*) \end{array} \right]
\end{eqnarray}
Note that  the rigidity matrix for $(G,p)$
can be expressed in block form
\begin{eqnarray}
\label{Eqn2}
R_q(G,p) = \left[ \begin{array} {cc}
X_0 & X_2(p^*) \end{array} \right]
\end{eqnarray}
Since $X_2(p')$ has maximal rank and $(G,p)$ is infinitesimally rigid we have \[\dim \ker X_2(p')\leq\dim \ker X_2(p^*)=0\] 
We conclude that  $u_2=u_3=\cdots=u_n=0$
and so $(G',p')$ is infinitesimally rigid.
\endproof

\begin{lemma}[({\cite[Theorem 3.2]{NOP2}})]
\label{GraphMoves}
A simple graph $G=(V,E)$ is $(2,2)$-tight if and only if there exists a finite sequence of graphs
\[K_1\to G^{(1)}\to G^{(2)}\to\cdots\to G\]
such that each successive graph is obtained from the previous by either a Henneberg $1$-move, a Henneberg $2$-move, a vertex to $4$-cycle move  or a vertex-to-$K_4$ move.
\end{lemma}

\proof[of Theorem \ref{qNormThm}]
If  $(G,p)$ is minimally infinitesimally rigid then $|E|=\rank R_q(G,p)
=2|V|-2$.
Also, if $H=(V(H),E(H))$ is a  subgraph of $G$ then the rows of its rigidity matrix $R_q(H,p)$ are linearly independent.
Thus we have the count $|E(H)|=\rank R_q(H,p)\leq 2|V(H)|-2$ and so $G$ is $(2,2)$-tight.

Conversely, if $G$ is $(2,2)$-tight then by, Lemma \ref{GraphMoves}, $G$ can be constructed inductively by applying finitely many graph moves
$K_1\to G^{(1)}\to G^{(2)}\to \cdots \to G$.
By the previous lemmas there exist  realisations $p^{(j)}$ for each $G^{(j)}$ such that $(G^{(j)},p^{(j)})$ is infinitesimally rigid.
In particular, there exists such a realisation for $G$ and so all  regular realisations of $G$ are infinitesimally rigid. 
Let $G'=(V,E')$ be a graph obtained by removing a single edge from $G$.
Then $\rank R_q(G',p)\leq|E'|=2n-3$ and so $(G',p')$ is not infinitesimally rigid by Proposition \ref{RigidityMatrix}. 
We conclude that $(G,p)$ is minimally infinitesimally rigid.
\endproof

\begin{remark}
We conjecture that the graphs $G$ which are generically infinitesimally rigid for the non-Euclidean spaces $(\mathbb{R}^d, \|\cdot \|_q)$ with $1<q<\infty$ are the $(d,d)$-tight graphs. 
This result would be somewhat analogous to the characterisation
of generically rigid body-bar frameworks given by Tay \cite{Tay-bodybar} who characterises  such frameworks in all dimensions in terms of
$(k,k)$-tight multigraphs (where $k=d(d+1)/2).$ 
We also remark  that higher dimensional normed spaces, including normed spaces of  matrices and infinite-dimensional spaces provide further intriguing contexts for the analysis of generalised bar-joint frameworks. 
There seem to be two natural approaches for this, namely (i) the development of inductive techniques  as we have done here, which lends itself well to simple (non-multigraph) graph settings, and (ii) the development of a more matroidal approach via submodular functions and spanning tree decompositions in the manner of 
\cite{F&Z}\cite{JJ_2005}\cite{N-W_1964}\cite{Tay3}\cite{Tay1999}\cite{Whi5}.
\end{remark}

%%%%%%%%%%%%%%%%%%%%%%%%%%%%%%%%%%%%%%%%%%%%%%%%%%%%%%%%%%%%%%%%%
\section{Polytopic norms}
\label{PolytopicNorms}
In this section we use a graph colouring technique to characterise the minimally infinitesimally rigid bar-joint frameworks in $\mathbb{R}^2$ with respect to certain norms for which the unit ball is a convex polytope. 
To define such a norm choose a spanning set of distinct non-zero vectors $b_1,b_2,\ldots,b_s\in \mathbb{R}^d$ and let $\mathcal{P}$ be the polytope
\[\mathcal{P}=\bigcap_{k=1}^s \{a\in\mathbb{R}^d:|a\cdot b_k|\leq 1\}\]
The associated polytopic norm on $\mathbb{R}^d$ is 
$\|a\|_\mathcal{P} = \max_{1\leq k\leq s} |a\cdot b_k|$.
As special cases we obtain both the $\ell^1$ norm and the $\ell^\infty$ norm.

For each $a=(a_1,\ldots,a_d)\in\mathbb{R}^d$ we define $\kappa(a)=b_k$ if there is a unique index $k$ such that $\|a\|_\mathcal{P}=|a\cdot b_k|$  and $\kappa(a)=0\in\mathbb{R}^d$ otherwise. Let $(G,p)$ be a bar-joint framework in the normed space $(\mathbb{R}^d,\|\cdot\|_\mathcal{P})$.

\begin{definition}
The {\em rigidity matrix} $R_\mathcal{P}(G,p)$ is a $|E|\times nd$ matrix with rows indexed by the edges of $G$ and $nd$ many columns indexed by the coordinates of the vertex placements $p_1,\ldots,p_n$. The row entries which correspond to an edge $v_iv_j\in E$ are
\[
\left[\begin{array}{ccccccccccc}
 0 & \cdots& 0& \kappa(p_i-p_j) & 0 &\cdots & 0 & -\kappa(p_i-p_j) & 0 & \cdots & 0
 \end{array}\right]
 \]
 where  non-zero entries may only appear in the $p_i$ columns and the $p_j$ columns.
\end{definition}

As for the $\ell^q$ norms we say that a bar-joint framework in  $(\mathbb{R}^d,\|\cdot\|_\mathcal{P})$ is {\em regular} if the rank of the rigidity matrix $R_\mathcal{P}(G,p)$ is maximal over all realisations. Note that  the set of regular realisations of a bar-joint framework in $(\mathbb{R}^d,\|\cdot\|_\mathcal{P})$ form an open set but in contrast to the $\ell^q$ spaces it is no longer the case that the regular realisations are dense in $\mathbb{R}^{nd}$. This prevents us from adapting the rigidity preservation arguments of Section \ref{lqNorms}. Instead we will establish the preservation of a spanning tree property for infinitesimally rigid frameworks which is based on induced framework colourings.

Suppose $(G,p)$ satisfies the condition that for each edge $v_iv_j\in E$  there is a unique index $k\in\{1,2,\ldots,s\}$ such that $\|p_i-p_j\|_\mathcal{P}=|(p_{i}-p_{j})\cdot b_k|$. In this case we say the edge $v_iv_j$ has {\em framework colour} $k$ and that the framework $(G,p)$ is {\em well-positioned} in $(\mathbb{R}^d,\|\cdot\|_\mathcal{P})$. 

\begin{proposition}
\label{RigidityMatrix2}
Let $(G,p)$ be a well-positioned bar-joint framework in $(\mathbb{R}^d,\|\cdot\|_\mathcal{P})$.
Then 
\begin{enumerate}
\item
$u\in \mathbb{R}^{nd}$ is an infinitesimal flex for $(G,p)$  if and only if $R_\mathcal{P}(G,p)u=0$.
\item
$(G,p)$  is infinitesimally rigid if and only if 
$\rank R_\mathcal{P}(G,p) = dn-d$.
\end{enumerate}
\end{proposition}

\proof 
The reasoning is the same as the proof of Proposition \ref{RigidityMatrix} except that in this case for each edge $v_iv_j\in E$ we compute  
$\zeta_{ij}'(0)=(u_i-u_j)\cdot b_k$
where $k$ is the framework colour of the edge $v_iv_j$.
\endproof

For example the above proposition shows that any well-positioned realisation of $K_2$ and $K_3$ in  $(\mathbb{R}^2,\|\cdot\|_\mathcal{P})$  is infinitesimally flexible.

If $(G,p)$ is well-positioned then we call the induced edge colouring $\lambda_p:E\to\{1,2,\ldots,s\}$, $v_iv_j\to k$ where $\kappa(p_i-p_j)=b_k$ the {\em framework colouring} of $G$. We denote by $G_k$ the largest monochrome subgraph of $G$ for which each edge has framework colour $k$.

\begin{proposition}
\label{Polytope1}
 Let $(G,p)$ be a well-positioned bar-joint framework in $(\mathbb{R}^d,\|\cdot\|_\mathcal{P})$ such  that the framework colouring of $G$ only involves the colours $1,2,\ldots,m$ where $m\leq d$.
 If $(G,p)$ is infinitesimally rigid then each of the monochrome subgraphs $G_1,G_2,\ldots,G_m$ contains a spanning tree.
\end{proposition}

\proof
Let  $k\in\{1,2,\ldots,m\}$ and suppose that $G_k$ does not span the vertex set of $G$.
Then there exists a partition $V = V_1 \cup V_2$ for which there is
no $k$-coloured connecting  edge. 
Let $W$ be the $(m-1)$-dimensional subspace of $\mathbb{R}^d$ spanned by $\{b_i:i=1,\ldots,m, i\not=k\}$ and choose a non-zero vector $z\in W^\perp$.
Then define $u=(u_1,\ldots,u_n)\in\mathbb{R}^{nd}$ such that $u_j = 0$
if $v_j \in V_1$ and $u_j =z$ if $v_j \in V_2$. We claim that $u$ is a
non-trivial infinitesimal flex. Non-triviality is clear by Lemma \ref{RigidMotions}  since  $u$ has
both zero and non-zero components. Also the flex condition is satisfied for edges with both vertices in $V_1$ or in $V_2$. For
a connecting edge, $v_1v_2$ say, suppose the framework colour is $l (\not=k)$. 
We have $u_1=0$ and $u_2=z$ and for small enough values of $t$,
\[\|(p_1 + tu_1) - (p_2 + tu_2)\|_\mathcal{P}
=|((p_1-p_2)-tz)\cdot b_l|
=\|p_1-p_2\|_\mathcal{P}\]
This implies that $u$ is a non-trivial infinitesimal flex which contradicts our hypothesis.
We conclude that $G_k$ contains a spanning tree.
\endproof

\begin{proposition}
\label{Polytope2}
 Let $(G,p)$ be a well-positioned bar-joint framework in $(\mathbb{R}^d,\|\cdot\|_\mathcal{P})$  such  that the framework colouring of $G$  involves the colours $1,2,\ldots,m$ where $m\geq d$.
 If the monochrome subgraphs $G_1,G_2,\ldots,G_d$ each contain a spanning tree then $(G,p)$ is  infinitesimally rigid.
\end{proposition}

\proof
Suppose  $u=(u_1,\ldots,u_n)\in \mathbb{R}^{nd}$ is an infinitesimal flex for $(G,p)$.
Then $R_\mathcal{P}(G,p)u=0$ and from the rigidity matrix  we see that 
 if $G_1,G_2,\ldots,G_d$ each contain a spanning tree then $(u_i-u_j)\cdot b_k=0$  for each pair $i,j$ and each $k$.
We conclude that  the only infinitesimal flexes of $(G,p)$  are  trivial. 
\endproof

\begin{example}
\label{K4Ex}
Let $\mathcal{P}$ be a polytope in $\mathbb{R}^2$ defined by linearly independent vectors $b_1,b_2\in\mathbb{R}^2$.
Then every well-positioned and regular realisation of $K_4$ in $(\mathbb{R}^2,\|\cdot\|_\mathcal{P})$ is minimally  infinitesimally rigid.
To see this consider the realisation $p=(p_1,p_2,p_3,p_4)$ where  $p_1=0$, $p_2=b_1$, $p_3=b_1+(1-\epsilon)b_2$ and $p_4=(1+\epsilon)b_2$.
By Lemma \ref{Isometric} we may assume after a change of basis that $b_1=(1,0)$ and $b_2=(0,1)$.
 Note that if $\epsilon$ is sufficiently small and non-zero then the edges $w_1w_2$, $w_1w_3$ and $w_3w_4$  have framework colour $1$ and the edges $w_1w_4$, $w_2w_3$ and $w_2w_4$  have framework colour $2$. By the above proposition $(K_4,p)$ is infinitesimally rigid and minimally infinitesimally rigid by Proposition \ref{RigidityMatrix2}.
\end{example}

\begin{definition}
Let $G=(V,E)$ be a simple connected graph. A graph $G'=(V',E')$ is obtained from $G$ by a {\em vertex splitting move} by the following process:
\begin{enumerate}
\item Choose an edge $v_1v_2\in E$.
\item Adjoin a new vertex $v_0$ to $V$ so that $V'=V\cup\{v_0\}$ and adjoin the edges $v_0v_1$ and $v_0v_2$ to $E$.
\item Either leave any edge of the form $v_1u\in E$ unchanged or replace it with the edge $v_0u$. 
\end{enumerate}
\end{definition}

To characterise the minimally  infinitesimally rigid bar-joint frameworks in $(\mathbb{R}^2,\|\cdot\|_\mathcal{P})$ we will use the following result. 
The Henneberg $1$-move, Henneberg $2$-move and vertex-to-$K_4$ move were defined in Section \ref{lqNorms}.

\begin{lemma}[({\cite[Theorem 1.5]{N&O}})]
\label{GraphMoves2}
A simple graph $G=(V,E)$ is $(2,2)$-tight if and only if there exists a finite sequence of graphs
\[K_1\to G^{(1)}\to G^{(2)}\to\cdots\to G\]
such that each successive graph is obtained from the previous by either a Henneberg $1$-move, a Henneberg $2$-move, a vertex splitting move  or a vertex-to-$K_4$ move.
\end{lemma}

\begin{remark}
The above construction lemma differs from Lemma \ref{GraphMoves} in that vertex splitting moves are used in place of the vertex-to-$4$-cycle move.  
This induction scheme is somewhat more technical to establish, however, the local nature of the vertex splitting move readily allows for the construction of regular placements without recolouring.
\end{remark}

\begin{lemma}
\label{Colouring}
Let $\mathcal{P}$ be a polytope in $\mathbb{R}^2$ defined by linearly independent vectors $b_1,b_2\in\mathbb{R}^2$.
Let $(G,p)$ be a well-positioned bar-joint framework in $(\mathbb{R}^2,\|\cdot\|_\mathcal{P})$ such that the monochrome subgraphs $G_1$ and $G_2$ span the vertex set of $G$.
If $G'$ is obtained from $G$ by a  Henneberg $1$-move, a Henneberg $2$-move, a vertex splitting move or a vertex-to-$K_4$ move then
there exists a well-positioned realisation $(G',p')$ such that
 the monochrome subgraphs $G_1'$ and $G'_2$ each span the vertex set of $G'$.
\end{lemma}

\proof
If $G'$ is obtained from $G$ by a Henneberg $1$-move then we can choose $p'=(p_0,p_1,\ldots,p_n)$ where 
$p_0$ lies within a sufficiently small distance of  the intersection of the lines  $p_1+tb_1$ and $p_2+tb_2$.
The new edges $v_0v_1$ and $v_0v_2$ in $G'$ have framework colours $1$ and $2$ respectively and all other framework colours are unchanged.
Thus $G'$ is spanned by both monochrome subgraphs.

Suppose we apply a Henneberg $2$-move to $G$ based on the edge $v_1v_2\in E$ and the vertex $v_3\in V$.
Without loss of generality we can assume that  $v_1v_2$ has framework colour $1$.
Choose $p'=(p_0,p_1,\ldots,p_n)$ where $p_{0}$ lies within a sufficiently small distance of the intersection of the line through $v_1$ and $v_2$ with the line $v_3+tb_2$. The three new edges $v_0v_1$, $v_0v_2$ and $v_0v_3$ have framework colours $1$, $1$ and $2$ respectively.
All other framework colours are unchanged and so $G'$ is spanned by both monochrome subgraphs.

 Suppose that  a vertex splitting move is applied to $G$ based on the edge $v_1v_2\in E$ and the vertex $v_1\in V$.
 We assume without loss of generality that $v_1v_2$ has framework colour $1$.
 Choose $p'=(p_0,p_1,p_2,\ldots,p_n)$  where $p_0=p_1+\epsilon b_2$.
 For sufficiently small non-zero values of $\epsilon$ the new edges $v_0v_1$, $v_0v_2$  in $G'$ have framework colours $2$ and $1$ respectively.
 Each edge in $G$ of the form $v_1u$ is either left unchanged or is replaced with the edge $v_0u$. In either case the framework colouring is unchanged provided $\epsilon$ is sufficiently small. The framework colours are also unchanged for all remaining edges and so $G'$ is spanned by the monochrome subgraphs $G'_1$ and $G'_2$.

 If a vertex-to-$K_4$ move is applied to $G$ based on the vertex $v_1\in V$ then we can choose a placement for the new vertices of $G'$ in two steps: Firstly, we arrange for the placement of $K_4$ to have two monochrome spanning trees. We described such a placement in Example \ref{K4Ex}. Secondly, we scale this placement of $K_4$ by a factor $r$ so that all of its vertices are placed in a small neighbourhood of $p_1$. 
Choose $p'=(p_1',p_2',p_3',p_4',p_2,\ldots,p_n)$ where $p_1'=p_1$, $p_2'=p_1+rb_1$, $p_3'=p_2'+r(1-\epsilon)b_2$, $p_4'=p_1+r(1+\epsilon)b_2$.
By Lemma \ref{Isometric} we may assume that $b_1=(1,0)$ and $b_2=(0,1)$.
If $\epsilon$ is sufficiently small then the edges $w_1w_2$, $w_1w_3$ and $w_3w_4$  have framework colour $1$ and the edges $w_1w_4$, $w_2w_3$ and $w_2w_4$  have framework colour $2$. This completes the first step.
Every edge in $G$ of the form $uv_1$ is replaced with $uw_j$ for some $j$ and if
$r$ is sufficiently small then the framework colours for these replacement edges remain the same.
All other framework colours are unchanged and so $G'$ is spanned by both monochrome subgraphs.
\endproof

\begin{theorem}
\label{PolytopeNorm}
Let $\mathcal{P}$ be a polytope in $\mathbb{R}^2$ defined by a pair of linearly independent  vectors $b_1,b_2\in\mathbb{R}^2$ and let $(G,p)$ be a well-positioned  and regular bar-joint framework in $(\mathbb{R}^2,\|\cdot\|_\mathcal{P})$.
 The following statements are equivalent:
\begin{enumerate}
\item  $(G,p)$ is minimally infinitesimally rigid.
\item The monochrome subgraphs $G_1$ and $G_2$ induced by the framework colouring of $(G,p)$ are edge-disjoint spanning trees.
\item $G$ is $(2,2)$-tight.
\end{enumerate}
\end{theorem}

\proof
(i)$\implies$(ii).
If $(G,p)$ is minimally infinitesimally rigid then by Proposition \ref{Polytope1} the monochrome subgraphs $G_1$ and $G_2$ each contain a spanning tree. 
If we remove a single edge from $G$ then the corresponding framework is no longer infinitesimally rigid and so by Proposition \ref{Polytope2} each monochrome subgraph $G_k$ is itself a spanning tree.

 (ii)$\implies$(iii). The top count is clear since $|E|=|E(G_1)|+|E(G_2)|=2n-2$. If $H$ is a subgraph of $G$ then either $|V(H)|\leq 3$ or by adjoining edges to $H$ we can build a minimally infinitesimally rigid framework $(H',p)$. In either case we have $|E(H)|\leq 2|V(H)|-2$.

(iii)$\implies$(i).
If $G$ is $(2,2)$-tight then there is a finite sequence of graph moves $K_1\to G^{(1)}\to G^{(2)}\to \cdots \to G$ as described in Lemma \ref{GraphMoves2}. 
By repeated application of Lemma \ref{Colouring} we can choose well-positioned realisations  $(G^{(j)},p^{(j)})$ such that $G^{(j)}$ is spanned by both monochrome subgraphs $G_1^{(j)}$ and $G_2^{(j)}$.
 In particular, such a realisation exists for $G$  and so by Proposition \ref{Polytope2} $G$ has a well-positioned and infinitesimally rigid realisation.
All well-positioned regular realisations of $G$ must also be infinitesimally rigid by Proposition \ref{RigidityMatrix2} and so in particular $(G,p)$ is infinitesimally rigid.
Let $G=(V,E')$ be a graph obtained by removing a single edge from $G$.
Then $\rank R_q(G',p)\leq|E'|=2n-3$ and so by Proposition \ref{RigidityMatrix2} $(G',p)$ is not infinitesimally rigid. 
We conclude that $(G,p)$ is minimally infinitesimally rigid.
\endproof

%%%%%%%%%%%%%%%%%%%%%%%%%%%%%%%%%%%%%%%%%%%%%%%%%%%%%%%%%%%%%%%
\section{Further directions}
A general context for bar-joint frameworks is provided by a metric space $(\mathcal{M}, d)$ in which $\mathcal{M}$ is an embedded  manifold in $\mathbb{R}^n$ and where the metric is obtained by restriction of a general norm $\|\cdot \|_C$. The case $n=3$ and $\|\cdot \|_C= \|\cdot \|_2$ has been considered in \cite{NOP}\cite{NOP2} leading to  characterisations of infinitesimal rigidity for the circular cylinder, with two independent infinitesimal motions, and for surfaces of revolution, including the cone, the torus and elliptical cylinders, each of which has a single independent infinitesimal motion.
Adopting a non-Euclidean norm such as $\|\cdot \|_q$ generally reduces the number of infinitesimal motions of the surface. In particular the $\|\cdot \|_q$-cylinder 
$\mathcal{M}_q=\{(x,y,z): \|(x,y,0) \|_q = 1\}$ has a $1$-dimensional space of infinitesimal rigid motions.
It is in fact possible to obtain a characterisation of infinitesimal rigidity in this setting
in terms of $(2,1)$-tight graphs by arguing exactly as we have done here and employing the inductive characterisation in \cite{NOP2} of $(2,1)$-tight graphs - the base graph is $K_5\backslash e$ and with the four graph moves above there is an additional elementary graph move in which two graphs in the class are joined by a single edge.
However, for the unit sphere
$\mathcal{S}_q=\{(x,y,z): \|(x,y,z) \|_q = 1\}$ there are no non-trivial rigid motions. Establishing a rigidity theorem for frameworks which are vertex-supported by $\mathcal{S}_q$ along the lines above would require an inductive scheme for $(2,0)$-tight  graphs and this is not presently available.

%%%%%%%%%%%%%%%%%%%%%%%%%%%%%%%%%%%%%%%%%%%%%%%%%%%%%%%%%%%%%%%%%%

%%%%%%%%%%%%%%%%%%%%%%%%%%%%%%%%%%%%%%%%%%%%%%%%%%%%%%%%%%%%%%%%%%%
%\affiliationone{% in this example, two authors share an institution
%   D. Kitson and S.C. Power\\
%   Department of Mathematics and Statistics\\ Lancaster University\\ Lancaster LA1 4YF \\U.K. 
%	\email{d.kitson@lancaster.ac.uk\\ s.power@lancaster.ac.uk}}

\end{document}